\documentstyle{amsart}
\let\emptyset\undefined
\newsymbol\emptyset 203F   
\newsymbol\restr 1216      
\newsymbol\halfline 2048   
\newsymbol\unitint 2049    
\newsymbol\Jint 204A       
\newsymbol\KK 204B         
\newsymbol\LL 204C         
\newsymbol\Miod 204D       
\newsymbol\rationals 2051  
\newsymbol\reals 2052      
\def\Hstar{\halfline\mkern1mu{}^*}
\def\Mstar{\Miod^*}

\def\piH{\pi_\halfline}
\def\Lstar{\LL\mkern-5mu{}^*}
\def\Isubn{\unitint_n}
\def\Jintn{\Jint_n}
\def\Isubu{\unitint_u}
\def\Jintu{\Jint_u}
\def\Isubv{\unitint_v}
\def\Nstar{\omega^*}
\mathchardef\from="303A          
\def\fin{\mathrm{fin}}  
\def\cont{\frak{c}}
\mathchardef\calB="0242
\mathchardef\calC="0243
\mathchardef\mathhyphen "002D
\def\ulim{u\mathhyphen\lim}
\def\cl{\mathop{\mathrm{cl}}\nolimits}
\def\seq #1_#2{\langle #1_#2\rangle_#2}
\def\<#1,#2>{\langle#1,#2\rangle}
\def\CH{\mathrm{CH}}
\let\phi=\varphi
\def\citename#1{{\normalshape\sc#1}} 
\def\citenames#1#2{{\normalshape{\sc#1} and {\sc#2}}} 
\theorembodyfont{\normalshape\sl}
\newtheorem{theorem}{Theorem}[section]
\newtheorem{lemma}[theorem]{Lemma} 
\theorembodyfont{\normalshape\rm}
\newtheorem{remark}[theorem]{Remark}
\newtheorem{question}[theorem]{Question}
\def\proof{\par\removelastskip\smallskip
           \noindent \hbox to\parindent{$\Box$\hss}\ignorespaces}

\advance\leftmargini by -10pt
\advance\leftmarginii by -10pt
\makeatletter
\def\thebibliography{\section*\refname
                 \setbox0=\hbox{[1990]} 
                 \labelwidth=\wd0 \bibmargin=\labelwidth
                 \advance\bibmargin\leftoflabel
                 \advance\bibmargin\rightoflabel
                 \defaultfont\small
                 \frenchspacing\raggedright
                 \clubpenalty4000\relax \widowpenalty\clubpenalty
                 \sfcode`\.\@m}

\newdimen\bibmargin      
\newskip\betweenauthors     
\newdimen\safety           
\newdimen\leftoflabel    
\newdimen\rightoflabel   

\def\betweenauthorsskip{\par
                       \vskip \z@ plus \safety
                       \penalty-100
                       \vskip \z@ plus -\safety
                       \vskip\betweenauthors}
\def\bibitem[#1]#2 #3\newblock{
       {\def\protect##1{\string ##1\space}
        \immediate\write\@auxout{\string\bibcite{#2}{#1}}}%
       \def\thisauthor{#3}
       \ifx\thisauthor\lastauthor
           \par\penalty50 
       \else\betweenauthorsskip  
            \noindent{\sc#3}\par\nobreak%
       \fi
       \def\lastauthor{#3}
       \setbox0=\hbox{[#1]}%
       \ifdim\labelwidth<\wd0
       \else\setbox0=\hbox to \labelwidth{[#1]\hfil}%
       \fi
       \hangindent\bibmargin\noindent
       \hbox{\hskip\leftoflabel\box0\hskip\rightoflabel}%
       \ignorespaces}
\def\lastauthor{John Doe}
\makeatother

\begin{document} 

\title{\v{C}ech-Stone remainders\\of spaces that look like $[0,\infty)$}

\author{Alan Dow}
\address{Department of Mathematics\\York University\\4700 Keele Street\\
         North York, Ontario\\Canada M3J 1P3}
\email{dowa@@nexus.yorku.ca} 

\author{Klaas Pieter Hart}
\address{Faculty of Technical Mathematics and Informatics\\ TU Delft\\
         Postbus 5031\\2600~GA~~Delft\\the Netherlands}
\email{wiawkph@@dutrun2.tudelft.nl} 

\subjclass{Primary 54D40;
           Secondary 54F15 03C20}

\keywords{\v{C}ech-Stone remainder, continuum, real~line, Wallman spaces,
          elementary equivalence, saturated models}

\maketitle 

\begin{abstract}
We show that many spaces that look like the half~line~$\halfline=[0,\infty)$
have, under~$\CH$, a \v{C}ech-Stone-remainder that is homeomorphic
to~$\Hstar$.                             
We also show that $\CH$ is equivalent to the statement that all standard
subcontinua of~$\Hstar$ are homeomorphic.

The proofs use Model-theoretic tools like reduced products
and elementary equivalence.
\end{abstract}

\section*{Introduction} 

The purpose of this note is to answer (partially) some natural questions
about the \v{C}ech-Stone remainder of the real~line or rather the remainder
of the space~$\halfline=[0,\infty)$ as the remainder of~$\reals$ is just a
sum of two copies of~$\Hstar$.

\medskip
Our first result says that under~$\CH$ the space~$\Hstar$ is, to a
certain extent, unique:
if $X$~is a space that looks a bit like~$\halfline$ then $X^*$ and $\Hstar$
are homeomorphic.
To `look a bit like~$\halfline$' the space~$X$ must be a connected ordered
space with a first element, without last element, of countable cofinality and
of weight at most~$\cont$.
The weight restriction is necessary, because if the weight of~$X$ is larger
than~$\cont$ then so is the weight of~$X^*$ and therefore $X^*$ cannot
be homeomorphic with~$\Hstar$.

As a consequence various familiar connected ordered spaces have a
\v{C}ech-Stone remainder that is homeomorphic to~$\Hstar$.
So the remainders of the lexicographic ordered square (minus the vertical
line on the right) and of any Suslin line are homeomorphic to~$\Hstar$.

\smallskip
The second result is concerned with the so-called standard subcontinua
of~$\Hstar$: take a discrete sequence~$\seq I_n$ of closed intervals
in~$\halfline$ and put, for any~$u\in\Nstar$,
$I_u=\bigcap_{U\in u}\cl\left(\bigcup_{n\in U} I_n\right)$.
Then $I_u$~is a standard subcontinuum of~$\Hstar$.

We show that $\CH$ is equivalent to the statement that all standard
subcontinua are homeomorphic.
This solves Problem~264 from \citenames{Hart}{van Mill}
\cite{HartvanMill90}.

\smallskip
Our final result shows that certain subcontinua of the standard subcontinua
are homeomorphic to~$\Hstar$; the precise statement is in
Section~\ref{sec.more.continua}, here suffice it to say that these continua
are natural candidates for being homeomorphic to~$\Hstar$.

\medskip
As may be expected we shall not directly construct homeomorphisms between
the spaces in question---it's too hard to take care of $2^{\omega_1}$ points
in $\omega_1$~steps---but we show that the spaces have isomorphic bases
for the closed sets (isomorphic as lattices).
That this works follows from the results of \citename{Wallman}
from~\cite{Wallman38}, to be described in Section~\ref{sec.prelim}
below.

A few words on how we show that the bases are isomorphic as lattices:
We implicitly and explicitly use a powerful result from Model Theory
which says that under quite general circumstances various structures
are isomorphic.
In each case the bases are identified as reduced products of families
of easily described lattices.
The factors of these products are pairwise elementary equivalent
and hence so are the products themselves.
Furthermore these products satisfy an certain saturation property.
The combination of elementary equivalence and this saturation property
implies that the lattices are isomorphic.
A more detailed explanation can be found in Section~\ref{sec.more.continua}.

\medskip
The paper is organized as follows.
Section~\ref{sec.prelim} contains some preliminary remarks.
In Section~\ref{sec.main} we prove the first result, the proof is
self-contained (i.e., requires no model theory).
In Section~\ref{sec.more.continua} we prove the results about the standard
subcontinua, here we appeal to standard fact from Model Theory to keep
the proofs pleasantly short.
The final Section~\ref{sec.special} deals with a special case of
Theorem~\ref{thm.Hstar.unique} that can be proved under weaker assumptions.

\section{Preliminaries}\label{sec.prelim}

\subsection{Sums of compact spaces}\label{subsec.sums}
We shall be dealing with sums of compact spaces a lot, so it's worthwhile
to fix some notation.
So let $X=\bigoplus_{n\in\omega}X_n$ be a topological sum of compact spaces;
we always take $\bigcup_{n\in\omega}\{n\}\times X_n$ as the underlying set
of the space.
The map $q\from X\to\omega$ defined by~$q(n,x)=n$ extends
to $\beta q\from\beta X\to\beta\omega$.
We shall always denote the fiber of~$u\in\Nstar$ under the map~$\beta q$
by~$X_u$.

\subsection{The half line}\label{subsec.H.and.M}
Our main objects of interest are the half line $\halfline=[0,\infty)$
and its \v{C}ech-Stone remainder~$\Hstar$.
The space~$\Hstar$ is a quotient of another space---that is somewhat easier
to handle---by a very simple map.

Indeed, consider the space $\Miod=\omega\times\unitint$---the sum of
$\omega$~many copies of the unit interval~$\unitint$.
The map $\piH\from\Miod\to\halfline$ defined by~$\piH(n,x)=n+x$ maps $\Miod$
onto~$\halfline$ and the map $\piH^*=\beta\piH\restr\Mstar$ maps $\Mstar$
onto~$\Hstar$.

A key point in our proof is to see what kind of identifications are made
by~$\piH^*$, so we take a better look at the components of~$\Mstar$.
Because $\Nstar$ is zero-dimensional and because $\Isubu$ is connected for
every~$u$ we know exactly what the components of~$\Mstar$ are: the
sets~$\Isubu$.

Furthermore, each $\Isubu$ has a natural top and bottom: we call the point
$0_u=\ulim\<n,0>$ the bottom~point and $1_u=\ulim\<n,1>$ the top~point.
The continuum~$\Isubu$ has many cut~points: for every sequence~$\seq x_n$
in~$(0,1)$ the point $x_u=\ulim\<n,x_n>$ is a cut~point of~$\Isubu$ and this
set of cut~points is dense.
It follows that $\Isubu$ is irreducible between $0_u$~and~$1_u$, which means 
that there is no proper subcontinuum of~$\Isubu$ that contains 
$0_u$~and~$1_u$.

We can put a preorder on~$\Isubu$: say $x\le_uy$ if{}f every subcontinuum
of~$\Isubu$ that contains $0_u$~and~$y$ also contains~$x$.
The layer of the point~$x$ is the set $\{y:y\le_ux$~and~$x\le_uy\}$.
This order is continuous in the sense that $\{y:y\le_ux\}$ is the closure
of~$\{y:y<_ux\}$.
We shall use this order in the proof of Theorem~\ref{thm.ctble.cof.layer}.

We turn back to the map~$\piH^*$; using standard properties of the
\v{C}ech-Stone compactification one can easily prove the next lemma
(if $u\in\omega^*$ then $u+1$ is the ultrafilter generated by
$\{U+1:U\in u\}$).

\begin{lemma}\label{lemma.Hstar.from.Mstar}
For every $u\in\omega^*$ the map~$\piH^*$ identifies the points
$1_u$~and~$0_{u+1}$ and these are the only identifications made.
\end{lemma}       

The continua $\Isubu$ govern most of the structure of~$\Hstar$;
they are known as the standard subcontinua of~$\Hstar$.
More information on~$\Hstar$ can be found in the survey
\citename{Hart} \cite{Hart92}.

\subsection{Wallman spaces}\label{subsec.Wallman}

As mentioned above we construct the homeomorphisms indirectly via
isomorphisms between certain lattices of closed sets of the spaces in
question.

This is justified by the results of \citename{Wallman}
from~\cite{Wallman38}; Wallman generalized the familiar Stone duality for 
Boolean algebras and zero-dimensional spaces to a duality for lattices and 
compact spaces.
We briefly describe this `Wallman duality'. 

If~$L$ is a lattice then a filter on~$L$ is a subset~$F$ such
that $0\not\in F$, if $x_1,x_2\in F$ then $x_1\wedge x_2\in F$ and
if $x_1\in F$ and $x_1\le x_2$ then $x_2\in F$.
An ultrafilter on~$L$ is just a maximal filter.
The set $X_L$ of ultrafilters on~$L$ is topologized by taking the family of
all sets of the form $x^+=\{F:x\in F\}$ with~$x\in L$ as a base for the
closed sets of~$X_L$.
The space $X_L$ is always compact, it is Hausdorff if{}f $L$~satisfies a certain
technical condition.

If $\calB$ is a base for the closed sets of a compact Hausdorff space
then $\calB$ satisfies this condition.
Thus, $X=X_\calB$ whenever $\calB$~is a (lattice) base for the closed sets 
of~$\calB$.
It is now easy to see that two compact Hausdorff spaces with  
isomorphic (lattice) bases for the closed sets are homeomorphic.

\section{Remainders of spaces that look like~$\halfline$}
\label{sec.main}

This section is devoted to a proof of the result mentioned in the
introduction, namely

\begin{theorem}[$\CH$]\label{thm.Hstar.unique}
Let $X$ be a connected ordered space with a first element, with no last
element, of countable cofinality and of weight~$\cont$.
Then $X^*$ and $\Hstar$ are homeomorphic.
\end{theorem}

We shall construct the homeomorphism indirectly, via spaces that are
mapped onto~$\Hstar$ and~$X^*$ respectively.

Remember from~\ref{subsec.H.and.M} that $\Hstar$ is the quotient
of~$\Mstar$ obtained by identifying $1_u$~and~$0_{u+1}$ for
every~$u\in\Nstar$ and that the map is called~$\pi_\halfline$.

We can construct a similar situation for~$X^*$: take a strictly increasing
and cofinal sequence~$\seq a_n$ in~$X$ with $a_0=\min X$.
For every~$n$ let $\Jintn=[a_n,a_{n+1}]$ and consider the
sum~$Y=\bigoplus_n\Jintn$.
The map $\pi\from Y\to X$ defined by $\pi(n,x)=x$ identifies $\<n,a_{n+1}>$
and~$\<n+1,a_{n+1}>$ for every~$n$.

As in the case for $\Hstar$ and $\Mstar$ the only identifications made
by~$\pi^*=\beta\pi\restr Y^*$ are of $\ulim_n\<n,a_{n+1}>$
and~$\ulim_n\<n+1,a_{n+1}>={u+1}\mathhyphen\lim_n\<n,a_n>$ for
every~$u\in\omega^*$.
In other words, for every $u\in\Nstar$ the top point of~$\Jintu$ is
identified with the bottom point of~$\Jint_{u+1}$.
We denote the top~point of~$\Jintu$ by~$t_u$ and the bottom~point by~$b_u$.

This gives rise to the following lemma.

\begin{lemma}\label{lemma.about.inducing}
If $h\from\Mstar\to Y^*$ is a homeomorphism that maps $\Isubu$ to $\Jintu$
and moreover maps $1_u$ to~$t_u$ for every~$u$ then $h$~induces a
homeomorphism from~$\Hstar$ onto~$X^*$.
\end{lemma}

\begin{proof}
The maps $\pi\circ h$ and $\pi_\halfline$ have exactly the same fibers.
Both are closed, being continuous between compact spaces, hence quotient
mappings.
Hence $\Hstar$ (the quotient of~$\Mstar$ by $\pi_\halfline$) and~$X^*$
(the quotient of~$\Mstar$ by $\pi\circ h$) are homeomorphic.
\end{proof}

Our efforts then will be directed towards constructing a homeomorphism
between $\Mstar$ and $Y^*$ that satisfies the assumptions of
Lemma~\ref{lemma.about.inducing}.

Rather than constructing a homeomorphism we shall construct two bases
$\calB$ and~$\calC$ for the closed sets of~$\Mstar$ and~$Y^*$
respectively and an isomorphism between them that will induce the desired
homeomorphism.

To construct~$\calB$ we consider the lattice generated by the closed
intervals in~$\unitint$.
It is a base for the closed sets of~$\unitint$.

We let $\LL_n$ be the corresponding lattice for~$\Isubn$.
The product lattice~$\LL=\prod_n\LL_n$ corresponds in a natural way to a
base for the closed sets of~$\Miod$.
The reduced product $\Lstar=\prod_n\LL_n/\fin$---obtained by identifying
$x$~and~$y$ whenever $\{n:x(n)\neq y(n)\}$ is finite---will then correspond in
a natural way to a base for the closed sets of~$\Mstar$.
This will be the base~$\calB$.

In a similar way we find~$\calC$: let $\KK_n$ be the lattice generated by
the closed intervals of~$\Jintn$ and consider $\KK=\prod_n\KK_n$ and the
reduced product $\KK^*=\prod_n\KK_n/\fin$.
The lattice corresponds to a base~$\calC$ for the closed sets of~$Y^*$.

Finding an isomorphism between $\Lstar$ and $\KK^*$ is the same thing as
finding a bijection~$\phi$ between $\LL$ and $\KK$ such that for all
$x,y\in\LL$ we have $x\le^*y$ if{}f $\phi(x)\le^*\phi(y)$, where
$x\le^*y$ means that $\{n:x_n\le y_n\}$ is cofinite.

To ensure that the induced homeomorphism maps $\Isubu$ to~$X_u$ for
every~$u$, it suffices to ensure that whenever $y=\phi(x)$ the sets
$\{n:x_n=\emptyset\}$ and $\{n:y_n=\emptyset\}$ as well as the
sets $\{n:x_n=\Isubn\}$ and $\{n:y_n=X_n\}$ differ by a finite set only.

Furthermore, to get $h(0_u)=b_u$ and $h(1_u)=t_u$ for every~$u$ we simply
map the closed set $b_\Miod=\bigl\{\<n,0>:n\in\omega\bigr\}$ to
$b_X=\bigl\{\<n,a_n>:n\in\omega\bigr\}$
and the set $t_\Miod=\bigl\{\<n,1>:n\in\omega\bigr\}$ to
$t_X=\bigl\{\<n,a_{n+1}>:n\in\omega\bigr\}$.
We leave it to the reader to check that this will indeed suffice.

\medskip
We shall construct a bijection~$\phi$ from $\LL$ to~$\KK$ that satisfies
the following conditions:
\begin{itemize}
\item[($\alpha$)] $\phi(b_\Miod)=b_X$ and $\phi(t_\Miod)=t_X$, and
\item[($\beta$)] for every $x$ and $y$ in~$\KK$ there is an $N\in\omega$
                 such that for every $n\ge N$ the sets of endpoints of
                 $\phi(x)(n)$ and $\phi(y)(n)$ have the same configuration
                 as the sets of endpoints of $x(n)$ and $y(n)$.
                 By this we mean the following.
                 \begin{enumerate}
                 \item The closed sets $x(n)$ and $\phi(x)(n)$ have the same
                       number of intervals and the families of intervals
                       are similar in that if the $i$th interval of~$x(n)$
                       consists of one point then so does the $i$th interval
                       of~$\phi(x)(n)$ and vice versa.
                       The same is demanded of $y(n)$ and $\phi(y)(n)$.
                 \item If $\{a_i:i<k\}$, $\{b_j:j<l\}$, $\{c_i:i<k\}$ and
                       $\{d_j:j<l\}$ are the sets of endpoints of~$x(n)$,
                       $y(n)$, $\phi(x)(n)$ and~$\phi(y)(n)$ respectively
                       (all sets in increasing order) then for all $i<k$
                       and~$j<l$ we have $a_i\mathrel{<,=,>}b_j$ if{}f
                       $c_i\mathrel{<,=,>}d_j$.
                 \end{enumerate}
\end{itemize}
Condition~($\alpha$) is one of the demands made at the outset; in
combination with~($\beta$) it ensures that for example the sets
$\{n:x(n)=\Isubn\}$ and $\{n:\phi(x)(n)=X_n\}$ differ by a finite set.

Condition~($\beta$) also readily implies that $x\le^*y$ if{}f
$\phi(x)\le^*\phi(y)$ for all $x,y\in\KK$.

\medskip
By $\CH$ we can construct $\phi$ in an induction of length~$\omega_1$;
but rather than setting up the whole bookkeeping apparatus we show how to
perform a typical inductive step.
So assume we have a bijection $\phi\from A\to B$ that satisfies ($\alpha$)
and~($\beta$), where $A$ and $B$ are countable subsets of~$\KK$ and~$\LL$
respectively with $t_\Miod,b_\Miod\in A$ and $t_X,b_X\in B$.

Let $\seq x_i$ be an enumeration of~$A$ and let $y_i=\phi(x_i)$ for all~$i$.
We show how to find $\phi(x)$ for an arbitrary $x\in\KK\setminus A$
(the task of finding $\phi^{-1}(y)$ for $y\in\LL\setminus B$ is essentially
the same).

First find an increasing sequence $\seq n_k$ of natural numbers such that
whenever $i,j<k$ and $n\ge n_k$ the endpoints of $x_i(n)$~and~$x_j(n)$
and those of $y_i(n)$~and~$y_j(n)$ are in the same configuration.
Using the fact that the intervals $\Isubn$ and $X_n$ are densely ordered
it is now an easy matter to find $y\in\LL$ such that
the endpoints of $x(n)$~and~$x_i(n)$ and those of $y(n)$~and~$y_i(n)$
have the same configuration whenever $n_k\le n<n_k+1$ and $i<k$.
We put $\phi(x)=y$ of course.

\medskip
This completes the proof of Theorem~\ref{thm.Hstar.unique}.

\begin{remark}
Lemma~\ref{lemma.about.inducing} brings up an interesting question.
If $h\from\Mstar\to Y^*$ were just any homeomorphism then it would have to
map components of~$\Mstar$ to components of~$Y^*$ and thus would induce
a map~$\phi$ from $\Nstar$ to~$\Nstar$ by $h[\Isubu]=\Jint_{\phi(u)}$.
It is readily seen that $\phi$~is an autohomeomorphism of~$\Nstar$:
Note that the set $C=\{u:h(0_u)=t_{\phi(u)}\}$ is clopen so that we may 
change~$h$ by first turning the $\Isubu$ with $u\in C$ upside-down.
But then $\phi$ merely mirrors the action of~$h$ on the set~$\{0_u:u\in\Nstar\}$
and hence it is an autohomeomorphism.

The problem is now to find an autohomeomorphism of $\Mstar$ that permutes 
the~$\Isubu$ in the same way as $\phi^{-1}$ permutes the points of~$\Nstar$
for then we could simply say: if $\Mstar$ and $Y^*$ are homeomorphic
then~$\Hstar$ and ~$X^*$ are homeomorphic.
We formulate this as an explicit question.
\end{remark}

\begin{question}
Is there for every autohomeomorphism $\phi$ of~$\Nstar$ an 
autohomeomorphism~$h$ of $\Mstar$ such that $h[\Isubu]=\unitint_{\phi(u)}$
for all~$u\in\Nstar$?
\end{question}

\section{More homeomorphic continua}\label{sec.more.continua}

The argument given in Section~\ref{sec.main} is actually a careful proof
of a special case of a general Model-Theoretic result.
We shall give a brief sketch of this result and then show how it may be used
to show that a few more continua of interest are homeomorphic.

The result says ``elementary equivalent and countably saturated models
of size~$\omega_1$ are isomorphic''.

Two models for a theory are said to be elementary equivalent if they satisfy
the same sentences, where a sentence is a formula without free variables.
This may be rephrased in a more algebraic way;
two models $A$ and $B$ are elementary equivalent if{}f the following holds:
if $\{x_1,\ldots,x_n\}\subseteq A$ and $\{y_1,\ldots,y_n\}\subseteq B$
are such that for every formula $\phi$ with $n$~free 
variables $\phi(x_1,\ldots,x_n)$ holds in~$A$
if{}f $\phi(y_1,\ldots,y_n)$ holds in~$B$ then for every formula $\psi$ for
which there is an~$x\in A$ such that $\psi(x_1,\ldots,x_n,x)$ holds there
is also a~$y\in B$ such that $\psi(y_1,\ldots,y_n,y)$ holds (and vice versa
of course).

By way of example consider dense linear orders with first and last points.
Any two such sets are elementary equivalent:
if $F=\{x_1,\ldots,x_n\}$ and $G=\{y_1,\ldots,y_n\}$ are as in the previous
paragraph then we simply know that $x_i\le x_j$ if{}f $y_i\le y_j$
and $x_i$~is the first (last) element if{}f $y_i$~is.
The conclusion will then be: for every~$x$ that is in a certain position
with respect to~$F$ then there is a~$y$ in the same position with respect
to~$G$.

\smallskip
A countably saturated model is one in which, loosely speaking, every
countable system of equations has a solution if{}f every finite subsystem of
it has a solution.

A countably saturated dense linear order is generally known as
an~$\eta_1$-set: if $A$ and $B$ are countable and $a<b$ for every~$a\in A$
and~$b\in B$ then there is an~$x$ such that $a<x<b$ for all~$a$~and~$b$.

The well-known theorem of \citename{Hausdorff} from \cite{Hausdorff14}
that under~$\CH$ any two $\eta_1$-sets of cardinality~$\cont$ are isomorphic
can now be seen as a special case of the general isomorphism theorem.

The `typical inductive step' from Section~\ref{sec.main} may be modified
to show that that the reduced product modulo the finite sets is countably
saturated: we were looking for an element of~$\Lstar$ that satisfied the
same equations as~$x$ and we used the fact that we could always satisfy
any finite number of these equations.

\smallskip
We refer to the book \citenames{Chang}{Keisler} \cite{ChangKeisler77}
for the necessary background on Model Theory.

\medskip
We shall now use this Model-Theoretic approach to show that many more
continua are homeomorphic, under~$\CH$.
As noted in the introduction, the first result solves Problem~264 from
\citenames{Hart}{van Mill} \cite{HartvanMill90}.

\begin{theorem}
The Continuum Hypothesis is equivalent to the statement that all standard
subcontinua of\/~$\Hstar$ are homeomorphic.
\end{theorem}

\begin{proof}
One direction was done by \citename{Dow} in \cite{Dow84d}: under $\lnot\CH$
there are $u$~and~$v$ for which $\Isubu$ and $\Isubv$ are {\em not\/}
homeomorphic.

\smallskip
For the other direction we note that we can obtain a base for the closed
sets of~$\Isubu$ simply by taking the ultraproduct~$\LL/u$.
This product is actually an ultrapower because the~$\LL_n$ are all the same.

The proof is finished by noting that $\LL/u$ and $\LL/v$ are elementary
equivalent (both are elementary equivalent to~$\LL_0$)
and countably saturated (\citenames{Chang}{Keisler} 
\cite[Theorem~6.1.1.]{ChangKeisler77});
by the general isomorphism theorem the ultrapowers are isomorphic.
\end{proof}

The next result shows that, again under $\CH$, all layers of countable
cofinality are homeomorphic.
Indeed the following, stronger, theorem is true.

\begin{theorem}[$\CH$]\label{thm.ctble.cof.layer}
Let $\seq a_n$ be an increasing sequence of cut~points in some~$\Isubu$
and let $L$ be the `supremum' layer for this sequence.
Then $L$ is homeomorphic to~$\Hstar$.
\end{theorem}

\begin{proof}
To begin we note that, because $\Isubu$ is an $F$-space, the closed interval
$[a_0,L]$ is the \v{C}ech-Stone compactification of the interval~$[a_0,L)$.

We now follow the proof of Theorem~\ref{thm.Hstar.unique}.

Form the intervals $\Jintn=[a_n,a_{n+1}]$ and the topological sum
$Y=\bigoplus_n\Jintn$.
The map $\pi\from Y\to[a_0,L)$ that identifies
$\<n,a_{n+1}>$~and~$\<n+1,a_{n+1}>$ for every~$n$ induces a map from $Y^*$
onto~$L$ (the restriction of~$\beta\pi$).
This map is of the same nature as~$\pi_\halfline^*$: it identifies the
top~point of~$\Jintu$ and the bottom~point of~$\Jint_{u+1}$ for
every~$u\in\omega^*$.

Our aim is to find a homeomorphism $h\from\Mstar\to Y^*$ that satisfies
the assumptions of Lemma~\ref{lemma.about.inducing}.
We shall do this, again, via an isomorphism between bases for the closed
sets of~$\Mstar$ and $Y^*$ respectively.

We shall use the lattice~$\Lstar$ as a base for~$\Mstar$ and we make a base
for~$Y^*$ as follows:
For each~$n$ the interval~$\Jintn$ is homeomorphic with~$\Isubu$ and hence
it has a base~$\KK_n$ for the closed sets that is elementary equivalent
to~$\LL_n$.
The reduced product~$\KK^*=\prod_n\KK_n$ is then a base for the closed
sets of~$Y^*$.

We may now copy the inductive construction of a bijection from~$\LL$
to~$\KK$ from the proof of Theorem~\ref{thm.Hstar.unique}.
The only difference is that we can no longer rely on the linear order
of~$\Jintn$ when we are constructing the images coordinatewise.
Instead we enumerate the countably many formulas from lattice theory
with parameters from~$A$ and use elementary equivalence to produce for
every~$x$ a~$y$ such that, as $n$~gets bigger, there are more and more
formulas that $x(n)$~and~$y(n)$ both satisfy or both do not satisfy
(for the $y$'s we replace the parameters from~$A$ with their images
under~$\phi$ of course).
\end{proof}

\begin{remark}
The lattice $\KK_n$ is {\em not\/} the lattice generated by the intervals
of~$\Jintn$; indeed, the Wallman space of the latter lattice is a linearly
ordered continuum and in fact the continuum that one gets by collapsing the
layers of~$\Jintn$ to points.
\end{remark}

\section{A special case}\label{sec.special}

Consider the long line~$L$ of length $\omega_1\times\omega$; that is, we take
the ordinal~$\omega_1\times\omega$ and stick an open unit interval between
$\alpha$~and~$\alpha+1$ for every~$\alpha<\omega_1\times\omega$.

Apparently Eric van~Douwen raised the question whether $L^*$ and $\Hstar$
could be homeomorphic.
Theorem~\ref{thm.Hstar.unique} implies that the answer is yes,
under~$\CH$.

A slight modification of the methods in Section~\ref{sec.main} will show
that the answer is even yes if~$\frak{d}=\omega_1$.

\begin{theorem}
If\/ $\frak{d}=\omega_1$ then $L^*$ and\/ $\Hstar$ are homeomorphic.
\end{theorem}

\begin{proof}
We shall find, of course, a homeomorphism $h\from\Mstar\to Y^*$ of the
familiar kind, where $Y=\omega\times[0,\omega_1]$ and $[0,\omega_1]$ denotes
the long segment of length~$\omega_1$.

Now, because $\frak{d}=\omega_1$, we may take a sequence
$\langle x_\alpha:\alpha<\omega_1\rangle$
of points in~$\unitint^\omega$ with the following properties:
\begin{enumerate}
\item For all $\alpha$ and all $n$ we have $0<x_\alpha(n)<1$.
\item If $\beta<\alpha$ then $x_\beta<^*x_\alpha$.
\item If $x$ is such that $x(n)<1$ for all $n$ then there is an~$\alpha$
      such that $x<^*x_\alpha$.
\end{enumerate}

It is then an easy matter to define a sequence
$\langle h_\alpha:\alpha<\omega_1\rangle$ of homeomorphisms, where
$$
h_\alpha\from\bigcup_n\{n\}\times[0,x_\alpha(n)]\to\omega\times[0,\alpha],
$$
such that
\begin{itemize}
\item $h_\alpha\bigl[\{n\}\times[0,x_\alpha(n)]\bigr]=\{n\}\times[0,\alpha]$
      for all $\alpha$~and~$n$, and
\item if $\beta<\alpha$ then $h_\alpha$ extends $h_\beta$ except on a finite
      number of vertical lines.
\end{itemize}
It is then straightforward to check that this sequence induces the desired
homeomorphism from $\Mstar$ onto~$Y^*$.
\end{proof}

\begin{question}
Is $\frak{d}=\omega_1$ equivalent to the statement that $\L^*$ and $\Hstar$
are homeomorphic?
\end{question}

We note that $\frak{d}=\omega_1$ if{}f $\Mstar$ and $Y^*$ are homeomorphic;
this is so because $\frak{d}=\omega_1$ if{}f the character of the set of
top~points of~$\Mstar$ is~$\omega_1$.
We have just seen that this implies that $\Mstar$ and~$Y^*$ are
homeomorphic;
on the other hand if $\Mstar$ and $Y^*$ are homeomorphic then clearly the
set of top~points of~$\Mstar$ has character~$\omega_1$.


\end{document}